\documentclass[12pt,oneside]{article}
\usepackage{mathtext}
\usepackage[T2A]{fontenc}
\usepackage[cp1251]{inputenc}
\usepackage[russian]{babel}
\usepackage[dvips]{graphicx,rotating}
\usepackage{amssymb}
\usepackage{amsthm}
\usepackage{bm}
\usepackage{latexsym}
\usepackage{titlesec}
\usepackage{amsmath}
\usepackage{amsfonts}
\usepackage{cite}
\usepackage{eufrak}
\usepackage{eucal}
\usepackage{color}
\usepackage[a4paper, left=30mm, right=20mm, bottom=20mm, head=20mm]{geometry}
\newcounter{example}[section]
\newcommand{\exam}{\addtocounter{example}{1}}

 \numberwithin{equation}{section}

\begin{document}
\begin{center}
\textbf{\large{EQUATIONS OF THE FIRST KIND AND THE INVERSION OF  SERIES OF RESOLVENTS OF A CLOSED OPERATOR }}
\end{center}

\begin{center}
\textbf{A. R. Mirotin}\\
amirotin@yandex.ru
\end{center}

\

Abstract. Let $A$ be a densely defined closed operator in a complex Banach space $X.$ Conditions for left invertibility of operators
 of the form
$\sum_{j=1}^\infty a_j (\alpha_j -A)^{-1}$
are given. Several examples are considered.

\

Key wards: closed operator, inverse operator, ill-posed problem, integral equation of a first kind, inverse problem, signal processing.

\


\

\begin{center}
\textbf{\large{УРАВНЕНИЯ ПЕРВОГО РОДА И ОБРАЩЕНИЕ РЯДОВ ИЗ
РЕЗОЛЬВЕНТ ЗАМКНУТОГО ОПЕРАТОРА}}
\end{center}

\begin{center}
\textbf{А.Р. Миротин}
\end{center}

Аннотация. Даны условия левой обратимости операторов
вида
$\sum_{j=1}^\infty a_j(\alpha_j -A)^{-1},$
где $A$  есть плотно определенный замкнутый оператор в комплексном банаховом пространстве  $X.$ Рассмотрены примеры.

\

\centerline{\textbf{1. Введение}}

\

Как известно, обратные задачи, возникающие в
различных разделах науки, приводят к разнообразным уравнениям
первого рода, решение которых представляет собой, как правило, некорректную задачу (см., например, \cite{IVT}, \cite{LRS}, \cite{Kaban}, \cite{AG}). Рассмотрим следующее интегральное уравнение первого
рода в пространстве $L^p(\mathbb{R})$ ($1\leq p<\infty$):
$$
\int_t^\infty k(s-t)x(s)ds = y(t).\eqno(1)
$$
 Пусть функция $k$ разлагается в
ряд экспонент
$$
k(t)=\sum_{j=1}^\infty a_je^{-\alpha_j t}, \ t\geq 0,
$$
причем  $\sum_{j=1}^\infty |a_j|<\infty,$ и
$\mathrm{Re}\alpha_j>0$ (по поводу условий разложимости в ряды экспонент см. \cite{Leon}, \cite[c. 199]{Sib}).
Тогда уравнение (1) приобретает вид
$$
\sum_{j=1}^\infty a_j\int_t^\infty e^{-\alpha_j (s-t)} x(s)ds =
y(t).
$$
Но известно \cite[c. 644]{DS} (см. также \cite[c. 363 --
368]{GGK}), что оператор дифференцирования $Dx=x'$ в
$L^p(\mathbb{R})$ с областью определения
$$
\mathrm{dom}(D)=\{x: x \mbox{ абсолютно непрерывна на  конечных
интервалах, } x'\in L^p(\mathbb{R}) \}
$$
 замкнут,   плотно определен, имеет спектр   $\sigma(D)=\imath \mathbb{R}$, и при
$\mathrm{Re}\alpha>0$ его резольвента есть
$$
(\alpha-D)^{-1}x(t)=-\int_t^\infty e^{-\alpha (s-t)} x(s)ds.
$$
Значит, (1)  можно записать в виде
$$
\sum_{j=1}^\infty a_j(\alpha_j -D)^{-1}x=-y.\eqno(2)
$$

В качестве второго примера рассмотрим  интегральное уравнение
$$
\int_{-\infty}^\infty k_1(s-t)x(s)ds = y(t)\eqno(3)
$$
в пространстве $L^2(\mathbb{R}).$ Предположим, что функция $k_1$ --- четная и разлагается в
ряд экспонент
$$
k_1(t)=\sum_{j=1}^\infty b_je^{- \imath \beta_j |t|},\ t\in \mathbb{R},
$$
причем  $\sum_{j=1}^\infty |b_j|<\infty,$ и
$\mathrm{Im}\beta_j<0.$ Тогда наше уравнение приводится к виду
$$
\sum_{j=1}^\infty b_j\int_{-\infty}^\infty e^{- \imath \beta_j |s-t|} x(s)ds =
y(t).
$$
С помощью преобразования Фурье легко проверить, что дифференциальный оператор $H:=-d^2/dt^2$ в пространстве $L^2(\mathbb{R})$
с областью определения
$$
\mathrm{dom}(H)=\{x: x^{\prime\prime}\in L^2(\mathbb{R})  \mbox{ в смысле обобщенных функций } \}
$$
(свободный гамильтониан) замкнут,   плотно определен, $\sigma(H)=\mathbb{R}_+$, и при
$\mathrm{Im}\beta <0$ его резольвента есть
$$
(\beta^2-H)^{-1}x(t)=\frac{ \imath}{2\beta}\int_{- \infty}^\infty e^{-\imath\beta |s-t|} x(s)ds.
$$
Значит, (3) можно записать в виде
$$
\sum_{j=1}^\infty (- 2\imath b_j\beta_j)(\beta_j^2 -H)^{-1}x=y,
$$
и мы снова получили уравнение типа (2).

В разделе 3 (пример 3) будет рассмотрена обратная задача теории
обработки сигналов, также приводящаяся к  уравнению такого типа.

 Данная работа посвящена решению некоторых уравнений вида (2), в которых оператор дифференцирования заменен произвольным
замкнутым плотно определенным оператором $A$ в комплексном
банаховом пространстве $X.$ В частности, показано, что данная задача является при неограниченном $A$ некорректной и  предложен подход к ее регуляризации.

 Легко проверяемая формула
$$
a_1(\alpha_1-A)^{-1}+a_2(\alpha_2-A)^{-1}=
(a_1+a_2)(\alpha_1-A)^{-1}\left(\frac{a_1\alpha_2+a_2\alpha_1}{a_1+a_2}-A\right)(\alpha_2-A)^{-1}
$$
показывает, что при условии $a_1,a_2> 0$ линейная комбинация двух
значений резольвенты  оператора $A$ обратима слева тогда и только тогда, когда
число $(a_1\alpha_2+a_2\alpha_1)/(a_1+a_2)$ не принадлежит
его точечному спектру $\sigma_p(A).$ Отсюда следует, что
всевозможные линейные комбинации с положительными коэффициентами
двух значений резольвенты $(\alpha_1-A)^{-1}$ и
$(\alpha_2-A)^{-1}$ оператора $A$ обратимы слева тогда и только
тогда, когда $\sigma_p(A)$ не пересекается с отрезком с концами $\alpha_1$ и
$\alpha_2,$  т.~ е. с выпуклой
оболочкой ${\rm conv}\{\alpha_1, \alpha_2\}$ множества
$\{\alpha_1,\alpha_2\}.$

Ниже мы установим аналогичное утверждение для  функций от операторов,
когда эти функции представимы в виде сумм рядов типа Вольфа-Данжуа (относительно
последних см. \cite{Sib}, \cite{L} и приведенную там библиографию) с
неотрицательными коэффициентами и, в частности, для рациональных
функций с простыми полюсами и положительными вычетами. Последний
случай для вещественных полюсов
 рассматривался в \cite{Obrresolvent}. При этом применялось
 функциональное исчисление, построенное в \cite{IZV1} --- \cite{IZV3}.
 Отметим также работу \cite{Trudy}, где рассматривался континуальный аналог
 этой задачи.

\

\centerline{\textbf{2. Основные результаты}}

\

 \textbf{Определение 1.} Пусть $\{\alpha_j\}$ --- не более чем
счетное  подмножество $\mathbb{C},$ $A$ --- плотно
определенный замкнутый оператор в комплексном банаховом
пространстве $X.$  Пусть функция $f$ в
некоторой  окрестности спектра $\sigma(A)$ допускает разложение в
ряда типа Вольфа-Данжуа
$$
f(z)=\sum_{j=1}^\infty \frac{a_j}{\alpha_j -z},\ a_j\in \mathbb{C} \eqno(4)
$$
(в случае
конечного множества $\{\alpha_j\}$ будем считать, что все члены
ряда начиная с некоторого номера равны нулю), причем
$$
\sum_{j=1}^\infty \frac{|a_j|}{\mathrm{dist}(\alpha_j,\sigma(A))}<\infty.
$$
 Тогда мы положим
$$
f(A):=\sum_{j=1}^\infty a_j(\alpha_j -A)^{-1}. \eqno(5)
$$

\

Если  $A$ --- плотно определенный замкнутый оператор в комплексном
банаховом пространстве $X,$ то через  $\cal{F}(A)$ будет
обозначаться пространство функций, голоморфных в некоторой (своей,
для каждой функции) окрестности множества $\sigma(A)$ и в
бесконечности. Напомним, что в соответствии с голоморфным
функциональным исчислением Рисса-Данфорда неограниченных
операторов \cite{DS} для каждой функции $f\in \cal{F}(A)$
определен ограниченный оператор
$$
f(A)=f(\infty)I+\frac{1}{2\pi\imath}\int_{\Gamma}f(\lambda)(\lambda-A)^{-1}d\lambda,
$$
где $I$ --- единичный оператор в  $X,$  а $\Gamma$ есть состоящая
из конечного числа спрямляемых жордановых кривых положительно
ориентированная граница открытого множества $V$, содержащего
спектр оператора $A,$ причем функция $f$ голоморфна на замыкании  $\overline{V}$
множества $V.$

\ Для доказательства основного результата нам понадобятся две
леммы.

 \textbf{Лемма 1.} \textit{В условиях определения 1  ряд
(5) сходится по норме оператора, и это определение  согласовано с
голоморфным функциональным исчислением Рисса-Данфорда
неограниченных операторов.}

Доказательство. Сходимость ряда (5) по норме оператора следует из
оценки
$$
\|(\alpha_j -A)^{-1}\|\leq
\frac{1}{\mathrm{dist}(\alpha_j,\sigma(A))}.
$$

Для доказательства второго утверждения леммы предположим дополнительно, что $f\in
\cal{F}(A)$ (отсюда, в частности, следует, что множество $\{\alpha_j\}$ ограничено), и выберем такое открытое подмножество $V$ комплексной
плоскости, граница $\Gamma$ которого состоит из конечного числа
спрямляемых жордановых кривых, что
$$
\sigma(A)\subset V\subset\overline{V}\subset \mathbb{C}\setminus
\overline{\{\alpha_j\}},
$$
и $f$ голоморфна на множестве $\overline{V}$ (черта обозначает
замыкание). Пусть контур $\Gamma$ имеет положительную ориентацию
относительно (быть может, неограниченного) множества $V.$ Тогда с
учетом того, что $f(\infty)=0,$ получаем, что значение $f(A)$ в
смысле голоморфного функционального исчисления Рисса-Данфорда есть
$$
\frac{1}{2\pi\imath}\int_{\Gamma}f(\lambda)(\lambda-a)^{-1}d\lambda=
$$
$$
\sum_{j=1}^\infty
a_j\frac{1}{2\pi\imath}\int_{\Gamma}\frac{1}{\alpha_j-\lambda}(\lambda-A)^{-1}d\lambda=
\sum_{j=1}^\infty a_j (\alpha_j-A)^{-1},
$$
что и утверждалось. При почленном интегрировании ряда мы
воспользовались тем, что
$$
\left\|\int_{\Gamma}\frac{1}{\alpha_j-\lambda}(\lambda-A)^{-1}d\lambda\right\|=
2\pi\|(\alpha_j-A)^{-1}\|\leq \frac{2\pi}{\mathrm{dist}(\alpha_j,\sigma(A))}.
$$
Лемма доказана.

\

\textbf{Замечание 1. } Пусть  $\{\alpha_j\}$ ---
счетное ограниченное подмножество $\mathbb{C},$ а жорданова область   $G$  в $\mathbb{C}$ такова, что $\overline{G}\subset \mathbb{C}\setminus
\{\alpha_j\}.$ Тогда любая функция  $f,$ голоморфная на $\overline{G},$ разлагается в ряд вида (4), но такое
разложение,  вообще говоря, не единственно  (см., например, \cite[\S 6]{Sib}). Из леммы 1 следует, что значение $f(A)$ не зависит от выбора этого разложения, т.~е. определение 1 корректно.

\

Нам понадобится также
следующее утверждение о нулях рядов вида (4).

 \textbf{Лемма 2.} \textit{Если функция  $f$  допускает разложение (4), где $a_j\geq 0,$  причем не все $a_j$ равны $0,$  то все корни
 уравнения $f(z)=0$ принадлежат
 замкнутой выпуклой оболочке ${\overline{\rm
conv}}(\{\alpha_j\})$ множества $\{\alpha_j\}.$}

Доказательство. Если допустить противное, то найдётся прямая на
комплексной плоскости, разделяющая ${\overline{\rm
conv}}(\{\alpha_j\})$ и некоторый корень $z_0$ этого уравнения.
Следовательно, найдётся прямая, разделяющая ${\overline{\rm
conv}}(\{\alpha_j-z_0\})$ и $0$. Совершая поворот  $z\mapsto
e^{\imath\varphi}z$ на подходящий угол $\varphi,$ получаем, что некоторая
прямая ${\rm Re}w=a,$ $a>0$ разделяет  ${\overline{\rm conv}}(\{
e^{\imath\varphi}(\alpha_j-z_0)\})$
 и $0$. Ясно, что
$$
\sum_{j=1}^\infty\frac{a_j}{w_j}=0,
$$
где $w_j=e^{\imath\varphi}(\alpha_j-z_0).$

Дробно-линейное преобразование $\zeta=1/w$ переводит прямую ${\rm
Re}w=a$ в окружность, проходящую через $0,$ внутри которой
расположен компакт (образ множества ${\overline{\rm conv}}\{w_j\}$
при этом преобразовании), содержащий все точки $\zeta_j:=1/w_j.$
Следовательно,
$$
\sum_{j=1}^\infty\frac{a_j}{w_j}=\sum_{j=1}^\infty a_j\zeta_j\ne 0
$$
($a_j\geq 0,$  причем не все $a_j$ равны $0$), и мы
получили противоречие.

\

\textbf{Теорема 1.} \textit{Пусть $\{\alpha_j\}$ --- не более чем
счетное ограниченное подмножество $\mathbb{C}$  и  $A$ --- замкнутый
плотно определённый оператор в комплексном банаховом пространстве
 $X,$ спектр которого не пересекается с замкнутой выпуклой оболочкой
 ${\overline{\rm
conv}}(\{\alpha_j\})$  множества $\{\alpha_j\}.$ Пусть функция $f$ разлагается в ряд (4), в котором $a_j\geq 0,$ $0<\sum_{j=1}^\infty a_j<\infty.$
Тогда левый обратный к оператору $f(A)$ существует, определен на $\mathrm{dom}(A)$  и имеет вид}
$$
f(A)^{-1}=\gamma +\beta A+h(A),
$$
\textit{где}
$$
\gamma=\frac{\sum_{j=1}^\infty
a_j\alpha_j}{\left(\sum_{j=1}^\infty a_j\right)^2}, \quad
\beta=-\frac{1}{\sum_{j=1}^\infty a_j},
$$
\textit{ функция $h(z):=1/f(z)-\gamma-\beta z$ принадлежит
$\mathcal{F}(A),$ и  $h(A)$ понимается в смысле голоморфного
функционального исчисления.}

Доказательство. В силу леммы 2 функция $g(z):=1/f(z)$  голоморфна
в проколотой окрестности бесконечности $\mathbb{C}\setminus{\overline{\rm
conv}}(\{\alpha_j\}).$ Рассмотрим следующие пределы:
$$
\beta=\lim\limits_{z\to\infty}\frac{g(z)}{z}=\lim\limits_{z\to
\infty}\frac{1}{zf(z)}=-\frac{1}{\sum_{j=1}^\infty a_j},
$$
$$
\gamma=\lim\limits_{z\to\infty}(g(z)-\beta
z)=\frac{1}{\sum_{j=1}^\infty
a_j}\lim\limits_{z\to\infty}\frac{\sum_{j=1}^\infty
 a_j\left( 1+\frac{z}{\alpha_j-z}\right)}{\sum_{j=1}^\infty \frac{a_j}{\alpha_j-z}}=
$$
$$
=\frac{1}{\sum_{j=1}^\infty
a_j}\lim\limits_{z\to\infty}\frac{\sum_{j=1}^\infty\frac{a_j\alpha_j}{\alpha_j-z}}
{\sum_{j=1}^\infty\frac{a_j}{\alpha_j-z}}=\frac{1}{\sum_{j=1}^\infty
a_j}\lim\limits_{z\to\infty}
\frac{\sum_{j=1}^\infty\frac{a_j\alpha_j}{\frac{\alpha_j}{z}-1}}{\sum_{j=1}^\infty\frac{a_j}{\frac{\alpha_j}{z}-1}}=
$$
$$
=\frac{\sum_{j=1}^\infty a_j\alpha_j}{\left(\sum_{j=1}^\infty
a_j\right)^2}
$$
(переход к пределу под знаком суммы ряда законен ввиду абсолютной
сходимости ряда  $\sum_{j=1}^\infty a_j\alpha_j$).

Функция $h(z):=g(z)-\gamma-\beta z$ принадлежит $\mathcal{F}(A)$
($h(\infty)=0$). Более того, оба слагаемых в правой части
очевидного равенства
$$
1=g(z)f(z)=(\gamma+\beta z)f(z)+h(z)f(z)
$$
принадлежат $\cal{F}(A)$ ($f$ имеет в бесконечности нуль первого
порядка). Следовательно, применяя функцию $(\gamma+\beta
z)f(z)+h(z)f(z)$ к оператору $A$ и воспользовавшись свойствами
голоморфного и полиномиального функциональных исчислений
\cite[VII.9]{DS}, будем иметь
$$
(\gamma +\beta A)f(A)+h(A)f(A)=I.
$$
 Таким образом, левый
обратный  $f(A)^{-1}$ существует и равен $g(A)=\gamma +\beta
A+h(A).$ А так как оператор $h(A)$ ограничен на $X$, то область
определения оператора  $f(A)^{-1}$ равна  $\mathrm{dom}(A)$, что и
завершает доказательство.

\

\textbf{Следствие 1.} \textit{Теорема 1 справедлива, в частности, для рациональных
функций с простыми полюсами $\alpha_j$ и положительными вычетами $a_j$ ($j=1,\dots,m$).
При этом функция $h$  имеет вид
$$
h(z)=\sum_{j=1}^m\sum_{k=1}^{m_j}\frac{c_{jk}}{(\alpha_j-z)^k}\quad (c_{jk}\in \mathbb{C}),
$$
и потому
$$
h(A)=\sum_{j=1}^m\sum_{k=1}^{m_j}c_{jk}(\alpha_j-A)^{-k}.
$$}

Доказательство.  Равенство для $h(z)$ есть разложение правильной дроби в сумму простых дробей.
Равенство для $h(A)$ теперь сразу следует из определения голоморфного
функционального исчисления Рисса-Данфорда.

\

\textbf{Следствие 2.} \textit{ Пусть выполнены условия  теоремы 1.}

1). \textit{Обратная
задача $f(A)x=y$ является корректной (по Адамару), если и только если  оператор
$A$ ограничен.}

2). \textit{ Обратная задача $f(A)x=y$ корректна по Тихонову на множестве $M\subset X$ (см., например, \cite{IVT}),  если и только если оператор $A$ непрерывен в относительной топологии множества $f(A)M.$}

Доказательство. Поскольку оператор $h(A)$ ограничен, оба утверждения следуют из того, что $\beta\ne 0.$

 \
\textbf{Замечание 2.} В силу лемм 1 и 2 существование  левого обратного к $f(A)$  в условиях теоремы 1
 следует также из \cite[с. 643, теорема 9]{DS}. Новизна теоремы 1
 состоит в том, что она дает удобное выражение для этого обратного оператора. В частности, из вида этого обратного следует, что  некорректность соответствующей обратной задачи полностью обусловлена членом $\beta A.$ Как показывает идущее ниже следствие 3, это позволяет в ряде случаев строить регуляризирующее семейство  для задачи $f(A)x=y.$

Напомним определение регуляризирующего  семейства  (см., например, \cite[c. 46]{LRS}, \cite{IVT}).

\textbf{Определение 2}. Пусть $K$ --- ограниченный оператор в $X.$ Семейство ограниченных операторов $R_\alpha:X\to X$ ($0<\alpha<\alpha_0$) называется
регуляризирующим для задачи $Kx=y,$ если  $R_\alpha Kx\to x\ (\alpha\to 0)$ при всех $x\in X.$

Например,  в случае гильбертова пространства для регуляризации по Тихонову \cite{TA} регуляризирующее семейство имеет вид ($K^*$ обозначает сопряженный оператор)
$$
R_\alpha =(\alpha m+K^*K)^{-1}K^*,
$$
где $m$ --- минимум функционала Тихонова
$$
J_\alpha x=\|Kx-y\|^2+\alpha \|x\|^2
$$
(см., например, \cite[c. 38]{Kir}).

\

\textbf{Следствие 3.} \textit{Пусть оператор  $A$ имеет ограниченный обратный  $K.$ Если $R_\alpha^0$ есть регуляризирующее семейство задачи $Kx=y,$ то
$$
R_\alpha=\gamma+\beta R_\alpha^0+h(A)
$$
есть регуляризирующее семейство задачи $f(A)x=y.$}

Доказательство. Так как $R_\alpha^0 A^{-1}x\to x\ (\alpha\to 0)$ при всех $x\in X,$ то $R_\alpha^0 y\to Ay\ (\alpha\to 0)$ при всех $y\in \mathrm{dom}(A).$ А так как $y:=f(A)x\in \mathrm{dom}(A),$ то
$$
\lim\limits_{\alpha\to 0}R_\alpha f(A)x=\gamma f(A)x+\beta\lim\limits_{\alpha\to 0}R_\alpha^0 f(A)x+h(A)f(A)x=
$$
$$
\gamma f(A)x+\beta Af(A)x+h(A)f(A)x=f(A)^{-1}f(A)x=x,
$$
что и требовалось доказать.

\

Как уже было отмечено во введении для случая, когда множество
$\{\alpha_j\}$  состоит из двух точек, если условие $\sigma(A)\cap
{\overline{\rm conv}}(\{\alpha_j\})=\emptyset$ предыдущей теоремы
не выполнено, левый обратный к $f(A)$ может не существовать.
Аналогичное утверждение верно для любого не более чем счетного множества $\{\alpha_j\}.$

 \textbf{Предложение 1.} \textit{Пусть $A$ ---
замкнутый плотно определённый оператор в комплексном банаховом
пространстве
 $X,$ точечный спектр которого  пересекается с  множеством
 ${\rm
conv}(\{\alpha_j\})\setminus\{\alpha_j\}.$ Тогда найдется
рациональная функция вида (4), для которой оператор $f(A)$ не
обратим слева.}

Доказательство. Пусть $\lambda\in \sigma_p(A)\cap {\rm
conv}(\{\alpha_j\}),$ $\lambda\notin \{\alpha_j\}.$ По
известной теореме Каратеодори найдутся такие неотрицательные числа
$k_1, k_2, k_3$ и числа $\alpha_{j_\nu}\in \{\alpha_j\}$, что
$\lambda=k_1\alpha_{j_1}+ k_2\alpha_{j_2}+ k_3\alpha_{j_3}$ \
 и $k_1+ k_2+ k_3=1.$ Пусть
$a_\nu:=k_\nu|\alpha_{j_\nu}-\lambda|^2.$ Тогда для функции
$$
f(z):=\frac{a_1}{\alpha_{j_1}-z}+\frac{a_2}{\alpha_{j_2}-z}+\frac{a_3}{\alpha_{j_3}-z}
$$
имеем
$$
\overline{f(\lambda)}=k_1(\alpha_{j_1}-\lambda)+k_2(\alpha_{j_2}-\lambda)+k_3(\alpha_{j_3}-\lambda)=0
$$
(черта обозначает комплексное сопряжение). Поэтому, если $x$ есть
собственный вектор оператора  $A,$ отвечающий собственному
значению $\lambda,$ то
$$
f(A)x=f(A)x-f(\lambda)x=\sum_{\nu=1}^3a_\nu((\alpha_{j_\nu}-A)^{-1}x-(\alpha_{j_\nu}-\lambda)^{-1}x)=
$$
$$
\sum_{\nu=1}^3a_\nu((\alpha_{j_\nu}-A)^{-1}(\alpha_{j_\nu}-\lambda)^{-1})(\lambda-A)x=0,
$$
что и завершает доказательство.

\

\centerline{\textbf{3. Примеры} }

\

Известные работы А. Ф. Леонтьева (см., например, \cite{Leon}) посвящены, в основном, рядам экспонент с неограниченными показателями. В связи с  примерами 1 и 2 отметим, что, как показано в  \cite[c. 199,  следствие 4]{Sib}, например,
 целые функции экспоненциального типа, меньшего единицы, разлагаются в абсолютно сходящиеся ряды экспонент \textit{с ограниченными показателями}.

\

\textbf{Пример 1.} Применим теорему 1 к уравнению
(1). При наложенных во введении на ядро $k$ условиях и в
предположении, что $a_j>0$ и множество $\{\alpha_j\}$ ограничено, получаем,
что уравнение (2), а вместе с ним и это уравнение
разрешимо тогда и только тогда, когда $y\in \mathrm{dom}(D)$, и
при таких $y$ имеет единственное решение
$$
x(t)=-\gamma y(t)-\beta y'(t)-h(D)y(t),
$$
где $\gamma,$ $\beta$ и $h(z)$ определены в теореме 1. В этом
примере  задача  является некорректной (по Адамару), причиной чему служит неограниченность оператора $D$ в пространстве  $L^p(\mathbb{R}).$ Следствие 3 позволяет найти регуляризирующее семейство задачи (1), если выбрано регуляризирующее семейство задачи
$D^{-1}x=y$ в  $L^p(\mathbb{R}).$

 \

\textbf{Пример 2.}  Аналогично, при наложенных во введении на ядро $k_1$ условиях и в
предположении, что $-\imath b_j\beta_j>0$ и множество $\{\beta_j\}$ ограничено, получаем,
что уравнение (3)
разрешимо тогда и только тогда, когда $y\in \mathrm{dom}(H)$, и
при таких $y$ имеет единственное решение
$$
x(t)=\gamma y(t)-\beta y^{\prime\prime}(t)+h(H)y(t),
$$
причем $\gamma,$ $\beta$ и $h(z)$ определяются по теореме 1, если положить $a_j= - 2\imath b_j\beta_j,$ $\alpha_j= \beta_j^2$ (применение преобразования Фурье, по-видимому, не позволяет получить этот результат).  В этом
примере  задача является некорректной (по Адамару) по причине  неограниченности оператора $H$ в пространстве  $L^2(\mathbb{R}).$ Как и в предыдущем примере, следствие 3 позволяет найти регуляризирующее семейство задачи (3), если выбрано регуляризирующее семейство задачи
$H^{-1}x=y$ в $L^2(\mathbb{R}).$

 \

\textbf{Пример 3.} Рассмотрим систему дискретного времени (рекурсивный
фильтр) конечного порядка,
описываемую разностным уравнением
$$
\sum_{k=0}^Nc_ky(n+k)=\sum_{l=1}^Nb_lx(n+l), \eqno(6)
$$
 связывающим  сигнал на входе
$x(n)$  с сигналом на  выходе $y(n)$ ($c_k, b_l\in \mathbb{C}$,
$n\in \mathbb{Z}$),
см., например,  \cite[c. 223]{S}. Предположим, что  характеристическое уравнение
$$
p(z):=\sum_{k=0}^Nc_kz^k=0
$$
имеет $N$ различных корней $z_j\in \mathbb{C}$ ($j=1,\dots,N$).
Пусть $q(z):=\sum_{l=1}^Nb_l z^l.$  Используя разложение на
простейшие дроби, получаем ($z\in \mathbb{C}\setminus\{z_j\}$)
$$
\frac{q(z)}{p(z)}=\sum_{j=1}^N\frac{a_j z}{z-z_j},\eqno(7)
$$
где
$$
a_j=\lim_{z\to z_j}\frac{q(z)}{zp(z)}(z-z_j).
$$
Если мы положим
$$
f(z):=\sum_{j=1}^N\frac{a_j }{z_j-z},
$$
то равенство (7) приобретает вид
$$
\sum_{l=1}^Nb_l z^l=\sum_{k=0}^Nc_k z^k (-zf(z)).\eqno(8)
$$
Пусть $X$ есть некоторое банахово пространство двусторонних
последовательностей (сигналов дискретного времени), в котором
действует оператор сдвига $T: x(n)\mapsto x(n+1).$ Для
применимости теоремы 1 нам достаточно предположить, что $a_j>0$ и
что $T$ плотно определен и замкнут в $X,$ $\sigma(T)\cap
\mathrm{conv}\{z_j\}=\emptyset$ (в теории обработки
сигналов $X$
--- это, как правило, пространство $\ell_p(\mathbb{Z}),$ его подпространство
$\ell_p(\mathbb{Z}_+)$
или другие подобные пространства, в которых $T$ ограничен и даже
изометричен). Поскольку функция $zf(z)$ принадлежит
$\mathcal{F}(T),$  из (8) следует, что при $x\in X$
$$
\sum_{l=1}^Nb_l T^lx=\sum_{k=0}^Nc_k T^k (-Tf(T))x. \eqno(9)
$$
Если мы положим $y:=(-Tf(T))x,$ то (9) превращается в (6), т.~е.
при наших предположениях мы нашли решение  $y$ уравнение (6) при
заданном $x.$ Заметим, что, например, в пространстве
$\ell_2(\mathbb{Z}),$ т.~е. в пространстве сигналов, энергия
которых конечна, это решение единственно. В самом деле, в этом
пространстве $\sigma(T)=\mathbb{T}$ (единичная окружность).
Поэтому полином $p(z)$ не имеет нулей на $\mathbb{T}.$ Применяя к
уравнению (6) преобразование Фурье
$$
x\mapsto \mathrm{F}x(e^{\imath
\theta}):=\sum_{n\in\mathbb{ Z}}x(n)e^{\imath n \theta}
$$
на группе $\mathbb{Z}$ (т.~е. переходя к  рядам Фурье), получаем
$p(z)\mathrm{F}y(z)=q(z)\mathrm{F}x(z)$ ($z\in \mathbb{T}$),
откуда и следует единственность. Далее для простоты  мы будем считать, что $X=\ell_2(\mathbb{Z}).$
 Тогда
обратная задача нахождения сигнала на входе по известному
сигналу на выходе равносильна решению уравнения
$$
f(T)x=-T^{-1}y.
$$
  В силу следствия 1 эта задача
имеет решение  при
всех $y\in \ell_2(\mathbb{Z}),$ и это решение имеет вид
$$
x=-\gamma T^{-1}y-\beta y-h(T)T^{-1}y,
$$
где $\gamma$  и  $\beta$  определены в теореме 1, в которой положено $\alpha_j=z_j,$
$h(z)=1/f(z)-\gamma-\beta z=-zp(z)/q(z)-\gamma-\beta z$ ---
рациональная функция. Рассматриваемая задача поставлена корректно (по Адамару).

\end{document}